\documentclass{gen-p-l}
\usepackage{latexsym}
\usepackage{amsfonts}
\usepackage{rawfonts}
\input{prepictex}
\input{pictex}
\input{postpictex}
\def\ZZ{\mathbb Z}
\def\RR{\mathbb R}

\def\CP{\mathbb C \mathbb P}
\def\bcp{{\mathbb C \mathbb P}}

\newtheorem{defn}{Definition}
\newtheorem{thm}{Theorem}
\newtheorem{prop}[thm]{Proposition}

\newtheorem{lem}[thm]{Lemma}

\def\Vol{\mbox{\rm Vol}}
\begin{document}
\sloppy
\title[Einstein Metrics, $4$-Manifolds 
\&   Differential Topology]{Einstein Metrics, Four-Manifolds,\\
and   Differential Topology}

\author{Claude LeBrun}
\address{Department of Mathematics\\
SUNY at Stony Brook\\
Stony Brook, NY 11794-3651}
\email{claude@@math.sunysb.edu}
\thanks{Supported 
in part by  NSF grant DMS-0305865.}

\date{}
\maketitle

A Riemannian metric $g$ on a smooth manifold $M$ is said to be {\em Einstein} 
if it has constant Ricci curvature, or in other words if 
 \begin{equation}
\label{al}
r = \lambda g, 
\end{equation}
where $r$ is the Ricci tensor of $g$ and $\lambda$ is some real constant  \cite{bes}. 
We still do not know if there are any obstructions to the
existence of Einstein metrics on high-dimensional manifolds, but it has now been 
known for three decades that not every $4$-manifold admits such
metrics  \cite{hit,tho}. Only recently, however, has it emerged 
that there are also obstructions to the existence of Einstein metrics which depend on the 
differentiable structure rather than just on the homotopy type of a $4$-manifold
\cite{lno,kot,lric,il1}. 
This article will attempt to give a concise  explanation of
this state of affairs. 

Our current understanding of the problem rests largely on certain 
curvature estimates which are deduced from the Seiberg-Witten equations. 
The strongest of these \cite{lric} was originally proved  indirectly, by invoking
 the solution to a generalized
version of the Yamabe problem. One  main purpose of 
the present article is to give (\S \ref{hope}) a new and simpler proof of this estimate,
using conformal rescaling only in order to introduce a 
 generalized form of the Seiberg-Witten equations. 
But this article will also attempt to clarify the nature of the resulting
obstructions, by systematically reformulating them in terms
of a new diffeomorphism invariant, called $\alpha (M)$, which 
is  introduced in \S \ref{monopoly}.

As   this article will make abundantly  clear,
Blaine Lawson's  work  on spin geometry and   scalar curvature   
has had a  deep and lasting impact on  my  own 
research.  On a more personal level,  Blaine  has  
also  been tremendous source of   inspiration   and 
encouragement throughout my many years at Stony Brook.  
I am lucky indeed to be able to call him 
 a friend and colleague, and it is a very  great pleasure for me    to  
be able to contribute an article to 
   this volume.

\section{Differential Geometry  on $4$-Manifolds}

The curvature and topology of $4$-manifolds are interrelated in 
a number of ways that have no adequate analogs in other dimensions. 
Many of these phenomena are intimately related to the fact that 
the bundle  $\Lambda^2$ of 2-forms over an
oriented Riemannian $4$-manifold $(M,g)$  invariantly  decomposes 
as the direct sum
\begin{equation} 
\Lambda^2 = \Lambda^+ \oplus \Lambda^- , 
\label{deco} 
\end{equation}
 of the eigenspaces of 
the Hodge star 
operator
$$\star: \Lambda^2 \to \Lambda^2.$$
The sections of $\Lambda^+$ are   characterized by 
$\star\varphi = \varphi$, and so are called {\em self-dual 
2-forms}, whereas the  
sections of $\Lambda^-$   satisfy  
$\star\varphi = -\varphi$, and are  called {\em anti-self-dual
2-forms}. 
Writing an arbitrary $2$-from uniquely as 
$$\varphi = \varphi^+ + \varphi^-,$$
where $\varphi^\pm \in \Lambda^\pm$, we then have
$$\varphi \wedge \varphi = \Big( |\varphi^+|^2 - |\varphi^-|^2\Big) d\mu_g , $$
where $d\mu_g$ denotes the metric volume form associated with the
fixed orientation.
The real  importance of  all this stems from the fact that  the curvature of
any connection on a vector bundle over  $M$ is a bundle-valued $2$-form, and 
(\ref{deco}) therefore always gives rise to a decomposition of any curvature
tensor into simpler pieces.

The Riemann curvature tensor of $g$ is a good 
case in point. 
By raising an index, 
 we can turn  the Riemannian curvature 
into  a self-adjoint linear map 
$${\mathcal R} : \Lambda^2 \to \Lambda^2$$
called the  {\em curvature operator}. 
Decomposing the 2-forms as in (\ref{deco}),
this linear endomorphism of $\Lambda^2$ can then
be decomposed into irreducible pieces 
\begin{equation}
\label{curv}
{\mathcal R}=
\left(
\mbox{
\begin{tabular}{c|c}
&\\
$W^++\frac{s}{12}$&$\stackrel{\circ}{r}$\\ &\\
\cline{1-2}&\\
$\stackrel{\circ}{r}$ & $W^-+\frac{s}{12}$\\&\\
\end{tabular}
} \right) . 
\end{equation}
Here $s$ denotes the {scalar curvature},
whereas 
$\stackrel{\circ}{r}=r-\frac{s}{4}g$ denotes the   
     trace-free part of the Ricci curvature.
The tensors $W_\pm$ are the trace-free pieces of the appropriate blocks,
and  are  respectively called the
{\em self-dual} and {\em anti-self-dual Weyl curvature} tensors. 
The corresponding pieces  of the 
Riemann  tensor enjoy the remarkable property 
of being {\em conformally invariant}, in the sense that they remain unaltered when the
metric is multiplied by an arbitrary smooth positive function. 
Notice that Einstein condition (\ref{al}) is satisfies iff  
${\mathcal R}$ commutes with the star operator; in more elementary terms, this amounts to the 
statement that a Riemannian $4$-manifold is Einstein iff the sectional
curvatures coincide for every orthogonal pair of $2$-planes $P,P^\perp\subset T_xM$,
$x\in M$.

Now let us  suppose that  $(M,g)$ is a {\em  compact} oriented
Riemannian  $4$-manifold.
The  Hodge theorem then tells us that every de Rham class
on $M$ has a unique harmonic representative, so that we have a canonical identification 
$$H^2(M,{\mathbb R}) =\{ \varphi \in \Gamma (\Lambda^2) ~|~
d\varphi = 0, ~ d\star \varphi =0 \} .$$
However,  the Hodge star operator $\star$ defines an involution of 
the right-hand side. 
We therefore have a direct-sum decomposition
\begin{equation}
	H^2(M, {\mathbb R}) = {\mathcal H}^+_{g}\oplus {\mathcal H}^-_{g},
	\label{harm}
\end{equation}
where
$${\mathcal H}^\pm_{g}= \{ \varphi \in \Gamma (\Lambda^\pm) ~|~
d\varphi = 0\} $$
are the spaces of self-dual and anti-self-dual harmonic forms.

The intersection form 
\begin{eqnarray*}
Q  :
H^{2}(M, {\mathbb R})\times H^{2}(M, {\mathbb R})	
 & \longrightarrow & ~~~~ {\mathbb R}  \\
	( ~ [\varphi ] ~ , ~ [\psi ] ~) ~~~~~ & 
	\mapsto  & \int_{M}\varphi \wedge \psi 
\end{eqnarray*}
becomes positive-definite when restricted to 
${\mathcal H}^+_{g}$, and  negative-definite when restricted to 
${\mathcal H}^-_{g}$. Moreover, these two subspaces are 
mutually orthogonal with respect to $Q$.
Indeed, combining an $L^{2}$-orthonormal  basis for 
${\mathcal H}^{+}_{g}$ with an $L^{2}$-orthonormal  basis for 
${\mathcal H}^{-}_{g}$ gives one a basis for $H^{2}({\mathbb R})$
in which the intersection form is
 represented by the diagonal matrix  	   
	    $$\left[ 
	  \begin{array}{rl}  \underbrace{
	 \begin{array}{ccc}
	   		 1 &  &	  \\
	   		  &	\ddots &   \\
	   		  &	 & 1
	   	  \end{array}}_{b_{+}(M)}
	   	   & 
	   		   \\ 
	  
	 {\scriptstyle b_{-}(M)}  \!
	    \left\{\begin{array}{r}
	    \\
	    \\
	    \\
	    \\
	    \end{array}
	    \right. \! \! \! \! \! \! \! \! \! \! \! \! \! \! \! 
	    &\begin{array}{ccc}
	   		 -1	&  &   \\
	   		  &	\ddots &   \\
	   		  &	 & -1
	   	  \end{array}
	    \end{array} 
	    \right] .
	    $$ 
The numbers
$b_{\pm} (M)= \dim {\mathcal H}^\pm_{g}$ are thus oriented 
homotopy invariants of $M$; namely, $b_{+}$ 
(respectively, $b_{-}$) is 
 the   maximal  dimension possible for   a linear  subspace
 of $H^{2}(M, {\mathbb  R})$ 
 on which the restriction of $Q$ 
is positive (respectively,  negative) definite. 

The assignment $g\mapsto {\mathcal H}_g^+$  gives one  an important  canonical map
$$\{\mbox{Riemannian metrics on }M\}\longrightarrow Gr^+_{b_+}[H^2(M,\RR)]$$
from the infinite-dimensional space of all metrics to 
the finite-dimensional Grassmannian of $b_+(M)$-dimensional 
subspaces of $H^2(M,\RR)$ on which the intersection form is 
positive definite. This map is called the
{\em period map} of $M$. It is obviously invariant under the
action of the identity component $\mbox{\it  Diff}_0(M)$
of the diffeomorphism group, and can also be shown to be  
 invariant under the action of the smooth functions 
$M\to \RR^+$ by 
conformal rescaling.  
A more subtle fact is that the
 period map is smooth, and has no critical points \cite{don}.

The difference $\tau (M)= b_{+}(M)-b_{-}(M)$ is called the 
{\em signature} of $M$.
Like the Euler characteristic
$\chi (M) = 2-2b_{1}(M) + b_{2}(M)$, it is the 
index of a geometric 
elliptic operator, and so may be expressed as a curvature integral. 
Indeed, to be explicit, one has 
  \begin{eqnarray*}
\chi (M)&=& \frac{1}{8\pi^2}\int_M \left[\frac{s^2}{24}+|W^+|^2
+ |W^-|^2  -\frac{|\stackrel{\circ}{r}|^2}{2}
\right]d\mu \\
\tau (M)&=& \frac{1}{12\pi^2}\int_M \Big[|W^+|^2
- |W^-|^2  \Big]d\mu 
\end{eqnarray*}
for absolutely any Riemannian metric $g$ on $M$.
In particular, we may combine these to obtain the Gauss-Bonnet-like  
formul\ae 
\begin{equation} 	\label{gb}
	(2\chi \pm 3\tau ) (M) 
=  \frac{1}{4\pi^2}\int_M \left[ \frac{s^2}{24}+2|W_\pm|^2
  -\frac{|\stackrel{\circ}{r}|^2}{2}
\right]d\mu.
\end{equation}
Notice, however,  that  
 the  integrand in this last expression is non-negative for any Einstein metric.
This gives us an  important   constraint, discovered 
independently by Thorpe \cite{tho} and Hitchin \cite{hit},  on the
topology of $4$-dimensional Einstein manifolds.

\begin{thm}[Hitchin-Thorpe Inequality]\label{ht}
If the smooth compact oriented 4-manifold
$M$ admits an Einstein metric $g$, then 
$$(2\chi + 3\tau )(M) \geq 0 ,$$
with equality iff $(M,g)$ is finitely covered by a flat $4$-torus $T^4$
or by $K3$  with a hyper-K\"ahler metric.  
Moreover, $(M,g)$ must also satisfy 
$$(2\chi - 3\tau )(M) \geq 0 ,$$
and this inequality is strict unless 
$(M,g)$ is finitely covered by a flat $4$-torus $T^4$
or by the orientation-reversed version of $K3$  with a hyper-K\"ahler metric.  
\end{thm}

Here a smooth Riemannian metric $g$ on an oriented $4$-manifold $M$
is called {\em hyper-K\"ahler} if the induced connection on $\Lambda^+$
is flat and trivial. Up to diffeomorphism, there is exactly one 
 simply connected compact $4$-manifold
 admits such a metric, called  $K3$, in honor
of Kummer, K\"ahler, and Kodaira.  This $4$-manifold is spin
(meaning that  its tangent bundle has  $w_2=0$), and has 
$b_+=3$, $b_-=19$. 
One  model of $K3$ is the quartic hypersurface
$$z_0^4 + z_1^4 + z_2^4 + z_3^4 =0$$
in $\CP_3$. The important point is that $K3$ is the 
underlying smooth oriented $4$-manifold of a simply
connected compact complex surface with $c_1=0$, because 
 a  truly 
 remarkable theorem of Kodaira asserts that all the complex surfaces satisfying these  conditions
are  deformation equivalent, and hence mutually diffeomorphic \cite{bpv}.

In order to put the Hitchin-Thorpe inequality to work,  let us consider the
following important way of constructing new manifolds.

\begin{defn}
Let $M_{1}$ and $M_{2}$ be smooth connected compact oriented $n$-manifolds. 
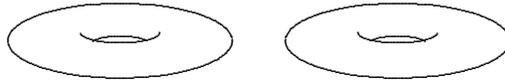
\begin{center}
\mbox{
\beginpicture
\setplotarea x from 0 to 290, y from 0 to 60
\ellipticalarc axes ratio 3:1  360 degrees from 140 40
center at 110 30
\ellipticalarc axes ratio 3:1  -360 degrees from 185 40
center at 215 30
\ellipticalarc axes ratio 4:1 -180 degrees from 125 33
center at 110 33
\ellipticalarc axes ratio 4:1 145 degrees from 120 30
center at 110 29
\ellipticalarc axes ratio 4:1 180 degrees from 200 33
center at 215 33
\ellipticalarc axes ratio 4:1 -145 degrees from 205 30
center at 215 29
\endpicture
}
\end{center}
Their connected sum $M_{1}\# M_{2}$ is then the smooth connected oriented $n$-manifold obtained 
by deleting a small ball from each manifold 
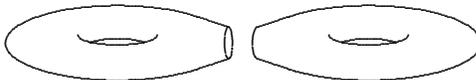
\begin{center}
\mbox{
\beginpicture
\setplotarea x from 0 to 290, y from 0 to 60
\ellipticalarc axes ratio 3:1  270 degrees from 145 40
center at 115 30
\ellipticalarc axes ratio 3:1  -270 degrees from 180 40
center at 210 30
\ellipticalarc axes ratio 4:1 -180 degrees from 130 33
center at 115 33
\ellipticalarc axes ratio 4:1 145 degrees from 125 30
center at 115 29
\ellipticalarc axes ratio 4:1 180 degrees from 195 33
center at 210 33
\ellipticalarc axes ratio 4:1 -145 degrees from 200 30
center at 210 29
\ellipticalarc axes ratio 1:4 360 degrees from 157 36
center at 157 30
\ellipticalarc axes ratio 1:3 180 degrees from 168 36
center at 168 30
{\setlinear 
\plot 145 40        157 36   /
\plot 145 20        157 24  /
\plot 180 40        168 36   /
\plot 180 20        168 24   /
}
\endpicture
}
\end{center}
and identifying the
resulting $S^{n-1}$ boundaries 
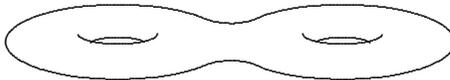
\begin{center}
\mbox{
\beginpicture
\setplotarea x from 0 to 290, y from 0 to 60
\ellipticalarc axes ratio 3:1  270 degrees from 150 40
center at 120 30
\ellipticalarc axes ratio 3:1  -270 degrees from 175 40
center at 205 30
\ellipticalarc axes ratio 4:1 -180 degrees from 135 33
center at 120 33
\ellipticalarc axes ratio 4:1 145 degrees from 130 30
center at 120 29
\ellipticalarc axes ratio 4:1 180 degrees from 190 33
center at 205 33
\ellipticalarc axes ratio 4:1 -145 degrees from 195 30
center at 205 29
{\setquadratic 
\plot 150 40    163 37    175 40   /
\plot 150 20    163 23    175 20  /
}
\endpicture
}
\end{center}
via a reflection. 
\end{defn}

If $M_1$ and $M_2$ are  simply connected $4$-manifolds, then $M= M_1\# M_2$
is also simply connected, and has  $b_\pm (M) = b_\pm (M_1) + b_\pm (M_2)$. 
Now let us use 
$\bcp_{2}$ denote the complex projective plane with its {\em standard} 
orientation, and $\overline{\bcp}_{2}$ denote the same smooth $4$-manifold
with the {\em opposite} orientation. 
Then the 
 iterated connected sum
 $$k \bcp_{2}\# \ell \overline{\bcp}_{2}
 =
\underbrace{\bcp_{2}\# \cdots \# \bcp_{2}}_{k } 
\# \underbrace{\overline{\bcp}_{2}\# \cdots \# \overline{\bcp}_{2}}_{\ell }$$
is a simply connected $4$-manifold with $b_+=k$ and $b_-=\ell$. In particular,
$$(2\chi + 3\tau ) (k \bcp_{2}\# \ell \overline{\bcp}_{2}) = 4 + 5k -\ell,$$
so this construction gives us  lots of simply connected $4$-manifolds which do not admit
Einstein metrics, by taking $\ell$ to be sufficiently large with respect to $k$. 

We  have now seen that  Euler characteristic $\chi$  and signature $\tau$
have a vital r\^ole to play in the theory of $4$-dimensional Einstein manifolds.  
However, knowing the  Euler characteristic  and 
signature  of a {\em simply connected} $4$-manifold is   equivalent to knowing 
 the invariants $b_\pm$, or in other words knowing
 the intersection form $Q$ up to isomorphism as a quadratic form
 over $\RR$.  Now a celebrated  result of Michael  Freedman \cite{freedman}
 asserts that  simply connected smooth $4$-manifolds are
 determined up to homeomorphism 
  by their intersection forms {over $\ZZ$}. However,  
 the integer coefficient version of the intersection form  
 certainly contains more information than just $b_\pm$; indeed, the 
 {\em parity} of the form (even or odd) determines 
 whether or not  the $4$-manifold is spin. 
 Indefinite quadratic forms over $\ZZ$ turn out to classified \cite{gost} by 
   parity and $b_\pm$. On the other hand,  
   a gauge-theoretic argument of Donaldson  \cite{donaldson}
shows that only the simplest definite forms 
can arise as intersection forms of $4$-manifolds. 
One thus obtains  the following remarkable classification result: 

\begin{thm}[Freedman]\label{fdmn}
Two smooth compact simply connected oriented $4$-manifolds are orientedly homeomorphic
if and only if 
\begin{itemize}
\item 
they have the same Euler characteristic $\chi$;
\item 
they have the same signature $\tau$; and 
\item both are spin, 
or both are 
non-spin.
\end{itemize} 
\end{thm}

As a consequence, any smooth compact simply connected  non-spin  $4$-manifold
is homeomorphic to a connected sum $k \bcp_{2}\# \ell \overline{\bcp}_{2}$.
For spin manifolds, the situation is a bit more unsettled, but the 
 connected sums 
$mK3\# n(S^2\times S^2)$ and their orientation-reversed versions, together with $S^4$,
at least exhaust all the  simply connected homeotypes 
for which  the additional constraint  
$\chi  \geq   \frac{11}{8}|\tau | +2$ is satisfied.
The so-called  $11/8$ conjecture  asserts that
this last inequality is in fact automatically satisfied,
or in other words  that the above list of spin homeotypes is complete. 
A  strong   partial result in this direction has been proved by Furuta 
\cite{fur108}.

For the purposes of  Riemannian  geometry, however, 
this beautiful classification is  somewhat aside from the point, since it is  concerned  with
 classification up to {\em homeomorphism},
not {\em diffeomorphism}. In order to do differentiable
geometry, we need 
a differentiable structure. However,  many of the above topological manifolds
turn out to 
have {\em infinitely many} inequivalent differentiable
structures, and  the  difference between these  smooth structures is often detectable 
by asking  questions about  the curvature 
of Riemannian metrics. 
 This article will attempt to   explain some of the 
 ramifications that  this interplay between  Riemannian geometry and 
differential topology  is now known to have 
for the existence and uniqueness of 
Einstein metrics on  $4$-manifolds.

\section{Examples of Einstein Manifolds}

Most of the currently available examples of 
compact $4$-dimensional Einstein manifolds are K\"ahler. 
Recall that a Riemannian  manifold
$(M,g)$ is called {\em K\"ahler} if it admits an 
almost-complex structure $J: TM\to TM$, $J^2= - 1$, which is 
invariant under parallel transport with respect to $g$.
Such an almost-complex structure is automatically integrable, and 
$(M,J)$ may therefore be viewed as a complex manifold. 

One may ask when a given compact complex manifold admits
a compatible K\"ahler metric which is also Einstein. For
Einstein metrics of negative Ricci curvature, the definitive
solution to this problem was found independently by Aubin \cite{aubin} and Yau \cite{yau}:

\begin{thm}[Aubin/Yau]\label{ay}
A compact complex manifold $(M^{2m},J)$ admits
a compatible K\"ahler-Einstein metric with
$s < 0$ iff its 
canonical line bundle $K=\Lambda^{m,0}$ is ample. When such a metric exists, it is 
unique, up to an overall multiplicative constant. 
\end{thm} 

Here a holomorphic line bundle $L\to M$  is called {\em ample} if it 
it has a positive power $L^{\otimes k}$ with enough holomorphic
sections to yield an embedding $M\hookrightarrow \CP_n$. 
As it turns out, 
a compact complex manifold $(M,J)$ of real dimension 
$4$ has ample canonical line bundle $K$ 
iff it is a minimal complex surface of general type without  $(-2)$-curves \cite{bpv}. 
These exist in great profusion.


Yau  \cite{yau} also gave  a definitive solution to  the analogous problem 
for Ricci-flat  metrics: 

\begin{thm}[Yau] \label{yau} 
A compact complex manifold $(M,J)$ admits 
a compatible K\"ahler-Einstein metric with 
$s =0$ iff $(M,J)$ admits a K\"ahler metric and 
$K^{\otimes \ell}$ is 
trivial for some positive integer $\ell$. 
When this happens, there is exactly one such metric
in each K\"ahler class.  
\end{thm} 

In real dimension $4$, there are exactly two diffeotypes of
compact K\"ahler manifolds for which $K$ is trivial, namely
$K3$ and $T^4$. The Ricci-flat K\"ahler metrics on these
manifolds are exactly the hyper-K\"ahler metrics alluded to in
our  previous discussion. Any other compact K\"ahler surface
with $K^\ell$ trivial for some positive $\ell$ is  the quotient of 
$K3$ or $T^4$ by the free action of a finite group of isometries
of one of these hyper-K\"ahler metrics.

The existence of K\"ahler-Einstein metrics is much more
delicate in the case of positive
 Ricci curvature. In real dimension $4$, however, 
 a complete solution to the problem was given by Tian \cite{tian}: 

\begin{thm}[Tian]\label{tian} 
A compact complex surface $(M^4,J)$ admits
a compatible K\"ahler-Einstein metric with
$s >  0$ iff
its anti-canonical line bundle $K^{-1}$ is ample and its 
Lie algebra of holomorphic vector fields is 
reductive. \end{thm} 

The  $4$-manifolds which carry Einstein metrics by virtue of 
this  last result are  
 ${\bcp}_2$, $S^2\times S^2$, and   
${\bcp}_2\#k\overline{\CP}_2$, $3\leq k \leq 8$.
  
 Until  quite recently, only sporadic examples of  
non-K\"ahler compact  $4$-dimensional Einstein  
manifolds were known. An interesting case in point is  the Page metric \cite{bes,page}
on $\CP_2\# \overline{\CP}_2$, which is beautiful, 
but has yet to lead to the construction of other compact examples. 
Of course, we have long known \cite{borel} that there is an infinite class of compact hyperbolic
manifolds ${\mathcal H}^4/\Gamma$, but it is easy to dismiss these
 constant-sectional-curvature spaces
as a bit boring in the small. 
 A recent construction of Michael Anderson \cite{mta}, however, 
 puts these hyperbolic manifolds  in a new, non-trivial context. Anderson's construction 
 begins with a complete, non-compact hyperbolic
 manifold of finite volume, replaces the cusps with 
 Schwarzschild-anti-deSitter  metrics, and then, under mild additional
 technical hypotheses,  perturbs the resulting
 metric  so as to make it Einstein. Of course, these  
 new Einstein manifolds still bear a
 family resemblance to their  hyperbolic 
 cousins, as they always have infinite fundamental group and  vanishing signature, 
  and always 
 contain large regions where  the sectional curvature is nearly 
 constant. Nevertheless, Anderson's construction does seem to represent the 
 first systematic method of constructing such a large class of compact
non-locally-symmetric  Einstein manifolds of general holonomy.

\section{The Seiberg-Witten Equations} 
\label{monopoly}

If $M$ is a smooth  oriented $4$-manifold,  $w_2(TM)\in H^2(M,\ZZ_2)$ is always in the 
image of the natural homomorphism $H^2 (M , \ZZ) \to H^2 (M , \ZZ_2)$ induced by 
$\ZZ \to \ZZ_2$. Consequently, 
we can always find
Hermitian line bundles $L\to M$ such that $c_{1}(L)\equiv w_{2}(TM)
\bmod 2$. For any such $L$, and for any Riemannian metric $g$ on $M$,
one can then find rank-$2$ Hermitian vector bundles ${\mathbb V}_\pm$
which formally satisfy
$${\mathbb V}_\pm= {\mathbb S}_{\pm}\otimes L^{1/2},$$
where ${\mathbb S}_{\pm}$ are the locally defined left- and
right-handed spinor bundles of $(M,g)$. 
Such a choice of   ${\mathbb V}_{\pm}$, up to isomorphism, 
is called a 
spin$^{c}$ structure $\mathfrak c$ on $M$.
Moreover, $\mathfrak c$ is completely determined by 
the first Chern class
$c_{1}(L)= c_{1}({\mathbb V}_{\pm})
\in H^{2}(M,\ZZ )$ if we assume that 
 $H_{1}(M,\ZZ)$ does not contain any elements of order $2$. 

Every 
unitary 
connection $A$ on $L$ induces a connection 
$$\nabla_A : \Gamma ({\mathbb V}_{+})\to \Gamma (\Lambda^1\otimes {\mathbb V}_{+}),$$
and composition of this with the natural {\em Clifford multiplication} homomorphism
$$\Lambda^1\otimes {\mathbb V}_{+}\to {\mathbb V}_{-}$$
gives one a 
spin$^c$ version 
$$D_{A}: \Gamma ({\mathbb V}_{+})\to \Gamma ({\mathbb V}_{-})$$
of the Dirac operator \cite{hitharm,lawmic}. This is an elliptic first-order differential operator, 
and in many respects closely resembles the usual Dirac operator of 
spin geometry. In particular,  one has  the Weitzenb\"ock formula 
\begin{equation}
\label{wtw} 
\langle \Phi , D_A^*D_A \Phi  \rangle = \frac{1}{2}\Delta |\Phi |^2 + |\nabla_A \Phi |^2 + 
\frac{s}{4} |\Phi |^2 + 2 \langle -iF_A^+ , \sigma (\Phi ) \rangle 
\end{equation}
for any  $\Phi\in \Gamma ({\mathbb V}_+)$, 
 where $F_A^+$ is the self-dual part of the 
curvature  of $A$,  and where 
$\sigma : {\mathbb V}_+ \to \Lambda^+$
is a natural real-quadratic map satisfying 
$$|\sigma (\Phi ) | = \frac{1}{2\sqrt{2}}|\Phi |^{2}.$$
This, of course, generalizes the Weitzenb\"ock formula 
used by  Lichnerowicz \cite{lic} to prove that metrics with $s > 0$ cannot exist
when $M$ is spin and $\tau (M) = -8\hat{A}(M)\neq 0$. However, one cannot hope to 
derive interesting  geometric information about the Riemannian
metric $g$ by just using (\ref{wtw}) for an arbitrary connection $A$,
since one would have no  control at all over the  
 $F_A^+$ term.
 Witten \cite{witten}, however,  had the brilliant insight that 
one could remedy this by  considering both $\Phi$ and $A$
as unknowns, subject to the  
{\em Seiberg-Witten equations}  
\begin{eqnarray} D_{A}\Phi &=&0\label{drc}\\
 F_{A}^+&=&i \sigma(\Phi) .\label{sd}\end{eqnarray}
These equations are non-linear, but they become an 
 elliptic first-order system 
once one imposes the `gauge-fixing' condition
$$d^* (A-A_0)=0$$
 to 
eliminate the natural action of the `gauge group' 
of  automorphisms of the Hermitian line bundle 
$L\to M$. 

Because the Seiberg-Witten equations are non-linear, one cannot 
use something like an index formula to  predict that they
must have solutions. Nonetheless,  there exist  spin$^{c}$ structures on 
many  $4$-manifolds  for which 
there is at least one solution 
for every metric $g$. This situation is conveniently described  by the 
following 
 terminology \cite{K}:
 
\begin{defn}
Let $M$ be a smooth compact oriented $4$-manifold
with $b_{+}\geq 2$. An element $a\in  H^{2}(M,\ZZ )/
\mbox{\rm torsion}$ will be called a {\bf monopole
class} of $M$ iff there exists a spin$^{c}$ structure
$\mathfrak c$ 
on $M$ with first Chern class 
$$c_{1}(L)\equiv a ~~~\bmod \mbox{\rm torsion}$$ which has the property
that  the corresponding  Seiberg-Witten 
equations (\ref{drc}--\ref{sd})
have a solution for every Riemannian  metric $g$ on $M$. 
\end{defn}

For example, if $(M,\omega )$ is a symplectic $4$-manifold with $b_+\geq 2$, and if 
$J$ is any almost-complex structure which is compatible with $\omega$, then 
$\pm c_1 (M,J)$ are both monopole classes \cite{taubes}. 
 Usually, one  detects the presence of  a monopole class
by  thinking of  the  moduli space of the
  Seiberg-Witten equations (that is, solutions 
  modulo gauge equivalence) as a cycle which represents an element of 
   the homology  
of a certain configuration space \cite{witten,taubes3,ozsz}. 
The resulting homology class is then metric-independent, and 
one may, for example,  then check that it is non-zero by analyzing the 
moduli space for some particular metric.  A    sophisticated recent refinment of this
idea, however,  instead 
detects the presence of a monopole class by means of an element of 
a stable  cohomotopy group  \cite{baufu,bauer2,il2}.

Because the Seiberg-Witten  equations imply the {Weitzenb\"ock formula}
 \begin{equation}
 0=	2\Delta |\Phi|^2 + 4|\nabla\Phi|^2 +s|\Phi|^2 + |\Phi|^4 ,	
 	\label{wnbk}
 \end{equation}
where $\nabla = \nabla_A$, 
one immediately sees that they admit no solution with $\Phi \not\equiv 0$
relative to a metric $g$ with  $s >0$. This allows one to prove, in particular, 
that there are lots of  
simply connected non-spin $4$-manifolds which do not admit 
positive-scalar-curvature metrics, in complete contrast to the situation in 
higher dimensions \cite{gvln}. 
But indeed, this Weitzenb\"ock formula makes an even more remarkable 
prediction concerning the behavior of the scalar curvature
\cite{witten,lpm}:

\begin{prop}\label{best}
Let $(M,g)$ be a smooth compact oriented Riemannian manifold,
let $\mathfrak c$ be a spin$^c$ structure on $M$, and let 
$c_1^+$ denote the self-dual part of the harmonic $2$-form 
representing the first Chern class $c_1(L)$ of $\mathfrak c$. 
If   there is a  solution of the Seiberg-Witten equations (\ref{drc}--\ref{sd}) 
on $M$ for $g$ and  ${\mathfrak c}$, then the scalar curvature 
$s_{g}$ of $g$ satisfies 
$$
\int_{M}s_{g}^{2}d\mu_{g} \geq 32\pi^{2} [c_{1}^{+}]^{2} .
$$
When $[c_1^+]\neq 0$, moreover,  equality can only occur if $g$
is a K\"ahler metric of constant scalar curvature. 
\end{prop}
\begin{proof}
Integrating (\ref{wnbk}), we have
$$0= \int [ 4|\nabla \Phi |^2 + s|\Phi|^2 + |\Phi|^4 ] d\mu , $$
and it follows that 
$$\int (-s) |\Phi|^2 d\mu \geq \int |\Phi|^4 d\mu .$$
Applying the Cauchy-Schwarz inequality to the
left-hand side therefore gives us 
$$
\left(\int s^2 d\mu  \right)^{1/2}\left(\int |\Phi |^4 d\mu  \right)^{1/2} \geq \int |\Phi |^4 d\mu , 
$$
so that 
$$
\int s^2 d\mu \geq \int |\Phi |^4 d\mu  = 8 \int |F_A^+|^2 d\mu ,
$$ 
and the inequality is strict unless  $\nabla \Phi \equiv 0$ and $s$ is constant. 
However, $F_A^+-2\pi c_1^+$  is an exact form plus a co-exact form, and 
so is $L^2$-orthogonal to the harmonic  forms. This gives us
the inequality 
$$\int |F_A^+|^2 d\mu \geq 4\pi^2 \int |c_1^+|^2 d\mu = 4\pi^2 \int c_1^+\wedge c_1^+,$$
and the last expression may be re-interpreted as the intersection pairing 
$[c_1^+]^2$ 
of 
the de Rham class of $c_1^+$ with itself. This gives us the desired
inequality
$$\int s^2 d\mu \geq 32\pi^2 [c_1^+]^2,$$
and, when the right-hand side is non-zero, equality can only happen if  
$g$ has special holonomy and constant scalar curvature.  
\end{proof}

One important application of this estimate is the following fundamental 
fact:

\begin{prop}
 Let $M$ be any smooth compact oriented $4$-manifold with 
 $b_+(M)\geq 2$. Then  ${\mathcal C}=\{ \mbox{monopole classes of }M\}$ 
  is a finite set. 
 \end{prop}
 \begin{proof}
 Let $g_1=g$ be any Riemannian metric on $M$, and let 
 ${\mathbf e}_1=[\omega_1]$ be the cohomology class of a harmonic self-dual
 form with respect to $g$, normalized so that ${\mathbf e}_1^2:=Q({\mathbf e}_1, {\mathbf e}_1)=1$.
 Because every metric is a regular point of the period map 
 \cite[Prop. 4.3.14]{don}, it follows that 
  $[\omega_1]$ has  an open 
  neighborhood in $H^2(M,\RR )$
in which every element can be represented by a self-dual harmonic form
relative to  some perturbation of $g$. However,  any open set in a finite-dimensional
vector space spans the entire 
space. Hence we can find a basis $\{ {\mathbf e}_j\}$ for $H^2(M,\RR )$, 
together with a collection of Riemannian metrics $\{ g_j ~|~j= 1, \ldots , b_2\}$, such that the 
$g_j$-harmonic $2$-form $\omega_j$ representing  the de Rham class ${\mathbf e}_j$ 
  is self-dual with respect  to $g_j$.  
  
  For convenience, let us  now normalize these basis
elements so that ${\mathbf e}_j^2=1$ for each $j$, and let $L_j: H^2(M,\RR ) \to \RR$
be the linear functionals $L_j(x)={\mathbf e}_j\cdot  x:= Q( {\mathbf e}_j,  x)$. 
Then 
Proposition 
\ref{best}, together with 
the Cauchy-Schwarz inequality, implies  that any monopole class 
$a\in H^2(M, \ZZ)/\mbox{torsion}$ satisfies 
$$|L_j (a)| = |{\mathbf e}_j \cdot a^+_{g_j}|\leq \sqrt{(a^+_{g_j})^2} \leq 
\left(\frac{1}{32\pi^2}\int_M s_{g_j}^2d\mu_{g_j}\right)^{1/2}=\kappa_j,$$
where the constant $\kappa_j$ is independent of $a$. This shows that 
 ${\mathcal C}\subset H^2(M,\RR )$
is contained in the $b_2(M)$-dimensional parallelepiped 
$$\left\{ x\in H^2 (M,\RR ) ~\Big|~ |L_j(x)|\leq \kappa_j ~\forall j= 1,\ldots , b_2(M)\right\},$$
which is  a compact set. 
Since ${\mathcal C}\subset H^2(M,\ZZ )/\mbox{torsion}$ is
also discrete, it follows  that $\mathcal C$ is finite. 
 \end{proof}

We now want to extract a numerical invariant from the set of 
monopole classes which captures those features of Seiberg-Witten
theory which are of the greatest relevance to problems in 
Riemannian geometry. To this end,  let us once again 
consider the open Grassmannian ${\mathbb G} = Gr^+_{b_+}[H^2(M, \RR )]$ of  all 
maximal linear subspaces ${\mathcal H}$  of the 
second cohomology for which the restriction $Q|_{\mathcal H}$ of the intersection pairing
  is positive definite. Each 
element ${\mathcal H}\in {\mathbb G}$ then determines 
an orthogonal  decomposition 
$$H^2(M, \RR ) = {\mathcal H} \oplus {\mathcal H}^\perp$$
with respect to $Q$. Let 
$${\mathfrak C}:={\mathcal C}-\{ 0\} \subset H^2(M,\ZZ )/\mbox{torsion}
\subset H^2(M,\RR)$$ denote the set of all 
the non-zero monopole classes $a\neq 0$ of $M$. Given a monopole class $a\in {\mathfrak C}$
and a positive subspace ${\mathcal H}\in {\mathbb G}$, we may then  define 
$a^+$ to be  the orthogonal projection  of $a$ 
into $\mathcal H$ with respect to $Q$. Using this, we now  define
an oriented-diffeomorphism invariant $\alpha (M)\in [0,\infty)$.
\begin{defn}
Let $M$ be a smooth compact oriented manifold with $b_+\geq 2$, and let
${\mathfrak C}\subset H^2(M, \ZZ ) /\mbox{torsion}$ be the set of
the non-zero monopole classes of $M$.  
If ${\mathfrak C}=\emptyset$, we declare that $\alpha (M)=0$.
Otherwise, we set
$$
\alpha (M) = \inf_{{\mathcal H}\in {\mathbb G}}\left[ \max_{a \in 
{\mathfrak C}} ~\sqrt{(a^+)^2} \right] ~ .
$$
\end{defn}
As a matter of notational convenience, let us  also set 
$\alpha^2(M):=[\alpha (M)]^2$. 

Needless to say, this is not the simplest invariant that one can
cook up using the Seiberg-Witten equations. However, 
it precisely captures those aspects of Seiberg-Witten theory
which are of the greatest relevance to many problems in Riemannian 
geometry. In particular, Proposition \ref{best} has the following
important consequence:

\begin{prop} Let $M$ be any smooth compact oriented $4$-manifold with $b_2\geq 2$. 
Then every Riemannian metric $g$ on $M$ satisfies
\begin{equation}
\label{yoho}
\int_M s_g^2 d\mu_g \geq 32\pi^2 \alpha^2 (M).
\end{equation}
Moreover, if $\alpha (M)\neq 0$, then 
 equality can hold only if $g$ is a K\"ahler metric of constant negative scalar curvature. 
\end{prop}

In particular, any metric, or sequence of metrics, on $M$ gives one an 
explicit upper bound for $\alpha (M)$. On the other hand, one can
obtain a lower bound for $\alpha (M)$ by replacing ${\mathfrak C}$ 
with any known set 
of non-zero monopole classes. Remarkably, the upper and lower bounds
obtained in this manner actually coincide for many $4$-manifolds, 
and in these circumstances one can therefore determine 
the invariant $\alpha$ 
even without necessarily knowing all the monopole classes on $M$. 
For example, the results of  \cite{lno,il2} allow one to read off the following: 

\begin{prop}
Let $M$ be the underlying smooth oriented $4$-manifold of any
compact complex surface $(M,J)$ with $b_+ > 1$. Let $X$ be the minimal model
of $(M,J)$. Then $$\alpha^2(M)=c_1^2 (X).$$ 
\end{prop}

\begin{prop}
Let $X_1$, $X_2$, $X_3$ be minimal, simply connected 
complex surfaces with $b_+\equiv 3\bmod 4$,
and let $M=X_1 \# X_2 \# X_3$. Then $$\alpha^2 (M) = c_1^2(X_1 ) + c_1^2 (X_2) 
+  c_1^2 (X_3).$$  
\end{prop}

Recall that a complex surface is said to be {\em minimal} if it is not obtained
from another complex surface by blowing up. Any complex surface $M$ can be
obtained from some minimal surface $X$, called its minimal model, by blowing
up a finite number of times  \cite{bpv}, and this minimal model $X$ is unique when $b_+> 1$.  
As smooth manifolds, one  then has
$$M\approx X \# k \overline{\CP}_2$$
for some integer $k \geq 0$, and the statement that $X$ is minimal amounts to saying
that this  value of $k$ is maximal. 

To see that all this is relevant to the study of Einstein metrics, notice that
(\ref{yoho}) can be rewritten as 
$$\frac{1}{4\pi^2}\int_M \frac{s^2}{24} d\mu \geq \frac{1}{3}\alpha^2 (M).$$
Comparing the left-hand side with our Gauss-Bonnet formula (\ref{gb})
for $(2\chi - 3\tau )(M)$, we thus  immediately obtain the following
improvement of one of the Hitchin-Thorpe inequalities:

\begin{thm}
If  a smooth compact oriented $4$-manifold $M$ 
with $b_+ > 1$  admits an  Einstein metric $g$, then 
$$
(2\chi - 3\tau ) (M) \geq \frac{1}{3} \alpha^2 (M).
$$
Moreover,
if $\alpha (M) \neq 0$,  equality occurs if and
only if $(M,g)$ is a compact quotient ${\mathbb C}{\mathcal H}_2/\Gamma$
of the complex hyperbolic plane, equipped with 
a constant multiple of its standard K\"ahler-Einstein metric.
\end{thm}
In particular, this tells us that the Einstein metric on a
complex-hyperbolic manifold ${\mathbb C}{\mathcal H}_2/\Gamma$
is {\em unique}, modulo diffeomorphisms and rescaling \cite{lno}.

Now one can obviously imitate this argument  by using $2\chi + 3\tau$ instead of
$2\chi - 3\tau$. However,  one would expect for the results obtained in this way  \cite{lno}
to be quite far from optimal, since the K\"ahler metrics which saturate
(\ref{yoho}) when $\alpha \neq 0$ always have $W_+\neq 0$. 
In the next
section, we will  remedy this, by 
replacing  Proposition
\ref{best} with an estimate which  involves  the self-dual Weyl curvature 
$W_+$ as well as the scalar curvature  $s$.

\section{The Main Curvature Estimate}
\label{hope}

The Dirac operator is conformally invariant, provided that the relevant
spinors are viewed as having appropriate conformal weights  \cite{lawmic,pr2}. 
Thus, if $\hat{\Phi}$ solves the spin$^c$ Dirac equation $\hat{D}_A\hat{\Phi}=0$
with respect to the conformally rescaled metric $f^{-2}g$, 
there is a corresponding solution $\Phi$ of the   equation $D_A\Phi=0$ 
with respect to $g$, with $|\Phi|_g= f^{-3/2}|\hat{\Phi}|_{\hat{g}}$ and 
$\sigma(\Phi) = f\hat{\sigma} (\hat{\Phi})$. If $(\hat{\Phi}, A)$
is a solution of the 
 Seiberg-Witten equations with respect to 
$f^{-2}g$, it thus follows that $(\Phi , A)$ solves the {\em rescaled Seiberg-Witten equations} 
\begin{eqnarray} D_{A}\Phi &=&0\label{rsdrc}\\
 F_{A}^+&=&if\sigma (\Phi) .\label{rssd}\end{eqnarray}
with respect to $g$. We will now show that directly studying 
equations (\ref{rsdrc}--\ref{rssd})   for an appropriate family of 
 choices of $f$ will yield an 
efficient avenue for deducing curvature estimates for the fixed metric $g$. 

Indeed,  plugging (\ref{rssd}) into the Weitzenb\"ock formula for (\ref{rsdrc}) gives us 
 \begin{eqnarray*}
0 & = & 2\Delta |\Phi |^{2} + 4 |\nabla_{A} \Phi |^{2} + s|\Phi |^{2}
 +8\langle -iF_A^+, \sigma (\Phi) \rangle  , \\
 & = & 2\Delta |\Phi |^{2} + 4 |\nabla_{A} \Phi |^{2} + s|\Phi |^{2}
 + f|\Phi |^{4} ,
\end{eqnarray*}
so that multiplying by $|\Phi |^2$ gives us 
$$
0=2|\Phi |^2\Delta |\Phi |^{2} + 4|\Phi |^2 |\nabla_{A} \Phi |^{2} + s|\Phi |^{4}
 + f|\Phi |^{6} .
$$
Integration therefore yields   
the  inequality 
\begin{equation}
\label{two}
0\geq \int_M\left[4 |\Phi |^{2}|\nabla_{A} \Phi |^{2} + s|\Phi |^{4}
 + f|\Phi |^{6}\right] d\mu .
\end{equation}
We will now use this to obtain an estimate involving $f$ and  curvature
of $g$.

 \begin{lem}\label{rescale}
 Let $(M,g)$ be a compact oriented Riemannian $4$-manifold, and let 
$f> 0$ be a smooth positive function on $M$. Suppose, for some
spin$^{c}$ structure with first Chern class  $c_{1}(L)$, that there is a 
 a solution of the Seiberg-Witten 
  equations with respect to the conformally related metric $f^{-2}g$. Let 
  $c_{1}^{+}$ denote the self-dual part of the harmonic 2-form representing $c_{1}(L)$,
relative to  the metric $g$.
 Then 
 \begin{equation}
 	\left(\int_M f^4d\mu_g \right)^{1/3}\left( \int_{M}\left|s_g+3w_g\right|^{3}f^{-2}d\mu_g 
 \right)^{2/3}\geq 72\pi^{2} [c_{1}^{+}]^{2},
 	\label{crux}
 \end{equation}
 where $s_g : M \to \RR$ denotes the scalar curvature of $g$ and 
   $w_g: M \to \RR$ is   the lowest eigenvalue of 
  $W_+: \Lambda^+\to \Lambda^+$ at each $x\in M$. 
 \end{lem}
 \begin{proof} 
 Any self-dual  2-form $\psi$ on any oriented 4-manifold satisfies 
 the 
Weitzenb\"ock formula \cite{bourg}
  \begin{equation}
\label{friend}
  (d+d^{*})^{2}\psi = \nabla^{*}\nabla \psi - 2W_{+}(\psi , 
\cdot ) + \frac{s}{3} \psi,
\end{equation}
where $W_{+}$ is the self-dual Weyl tensor. 
It follows that 
$$
  \int_{M} \left(|\nabla \psi |^{2}
-2W_{+}(\psi , \psi ) + \frac{s}{3}|\psi |^{2}\right) d\mu\geq 0, 	
$$
and hence that 
 $$
   \int_{M} |\nabla \psi |^{2}d\mu \geq 
\int_M\left(2w-\frac{s}{3}\right)|\psi |^{2} d\mu .
  $$ 
However, the  particular self-dual 2-form $\sigma (\Phi )$ satisfies 
\begin{eqnarray*}
	|\sigma (\Phi ) |^{2} & = & \frac{1}{8}|\Phi |^{4}  ,\\
	|\nabla \sigma (\Phi ) |^{2} & \leq  & \frac{1}{2} |\Phi |^{2}|\nabla \Phi |^{2} .
\end{eqnarray*}
 Setting $\psi = \sigma (\Phi )$,  
  we thus have
  $$
  \int_{M} 4|\Phi |^{2} |\nabla \Phi |^{2}
 ~d\mu \geq  \int_M\left(2w-\frac{s}{3}\right)|\Phi |^{4} d\mu .
   $$
 But   (\ref{two})
 tells us that 
 $$0\geq \int_M\left[4 |\Phi |^{2}|\nabla_{A} \Phi |^{2} + s|\Phi |^{4}
 + f|\Phi |^{6}\right] d\mu ,$$
 so we obtain 
 $$0\geq \int_M\left[ \left(\frac{2}{3}s+2w\right)|\Phi |^{4}
 + f|\Phi |^{6}\right] d\mu,$$
which we may rewrite as  
   $$
  \int_M\left[ -\frac{2}{3}\left(s+3w\right)f^{-2/3}\right]\left(f^{2/3}|\Phi |^{4}\right)
   d\mu \geq  \int_M f|\Phi |^{6}d\mu .
   $$
   Applying the H\"older inequality to the left-hand side now gives 
    $$
  \left[\int_M \left|\frac{2}{3}(s+3w)\right|^3f^{-2}d\mu\right]^{1/3}
  \left[ \int_M f|\Phi |^{6}d\mu \right]^{2/3}
   d\mu \geq  \int_M f|\Phi |^{6}d\mu ,
   $$
   and we therefore deduce that 
    $$
    \int_M \left|\frac{2}{3}(s+3w)\right|^3f^{-2}d\mu ~~\geq ~~ \int_M f|\Phi |^{6}d\mu .
    $$
    But the H\"older inequality also tells us that 
    $$
   \left(\int_M f^4d\mu\right)^{1/3}   \left(\int_M f|\Phi |^{6}d\mu\right)^{2/3} \geq 
   \int_M f^{4/3}(f^{2/3} |\Phi|^4) d\mu ~,    $$
   so it follows that 
  $$
   \left(\int_M f^4d\mu\right)^{1/3}   \left( \int_M \left|\frac{2}{3}(s+3w)\right|^3
   f^{-2}d\mu\right)^{2/3} \geq 
   \int_M f^2 |\Phi|^4 d\mu ~.    $$
   However, since $-iF_A^+ = f\sigma (\Phi )$, we also have 
   $$ \int_M f^2 |\Phi|^4 d\mu = 8 \int_M |F_A^+|^2 d\mu \geq 
  8 \int_M |2\pi c_1^+|^2 d\mu 
  =
   32\pi^2 [c_1^+]^2$$
   because $iF_A^+=2\pi c_1^+ + d\theta -  d^*(\star\theta)$ for some $1$-form $\theta$. Thus
   $$ \left(\int_M f^4d\mu\right)^{1/3}   \left( \int_M \left|\frac{2}{3}(s+3w)\right|^3f^{-2}d\mu\right)^{2/3} 
   \geq 32\pi^2 [c_1^+]^2.$$
   Multiplying both sides by $9/4$ now yields the promised inequality. 
 \end{proof}

   \begin{lem} \label{nearly} 
   Let $M$ be a smooth compact oriented $4$-manifold with $b_+(M) \geq 2$, and 
   let $a\in H^2(M,\ZZ )/\mbox{torsion}$ be a monopole class. Let $g$ be any
   Riemannian metric on $M$, and let $u>0$ be any smooth positive function on $M$ 
   for which 
   $$u \geq |s+3w|$$
   at each point $x\in M$. Then 
   $$\int u^2d\mu_g \geq 72\pi^2 (a^+)^2.$$ 
   \end{lem}
   \begin{proof}
   For any such function $u$, set $f=u^{1/2}$. and notice that  we have
   $$
   \left(\int u^2 d\mu\right)^{1/3} =  \left(\int f^4d\mu \right)^{1/3}
   $$
   and 
   $$
   \left( \int u^2 d\mu  \right)^{2/3}=  \left(\int u^3 f^{-2}d\mu  \right)^{2/3}\geq 
    \left(\int |s+3w|^3 f^{-2}d\mu  \right)^{2/3} .
   $$
   Setting $c_1(L)=a$ and invoking 
   Lemma \ref{rescale}, we  therefore have  
   $$\int_M u^2 d\mu_g \geq 
   \left(\int_M f^4d\mu_g \right)^{1/3}\left( \int_{M}\left|s+3w\right|_g^{3}f^{-2}d\mu_g 
 \right)^{2/3}\geq 72\pi^{2} (a^{+})^{2},
   $$
   as claimed.
   \end{proof}
   
   \begin{lem}\label{chiave}
   Let $M$ be a smooth compact oriented $4$-manifold with $b_+(M) \geq 2$, and 
   let $a\in H^2(M,\ZZ )/\mbox{torsion}$ be a monopole class. Then every 
   Riemannian metric $g$ on $M$ satisfies 
   $$\int \left(s+3w\right)^2 d\mu_g \geq 72\pi^2 (a^+)^2.$$ 
   \end{lem} 
  \begin{proof}
 Let $u_j$ be a sequence of smooth positive functions satisfying $u_j > |s+3w|$ and with 
 $u_j \to |s+3w|$ in the $C^0$ topology.  (Since the smooth functions are dense in $C^0$, 
 such a sequence may of course be constructed by 
 setting $u_j = \frac{1}{j} + v_j$, where $v_j$ is a smooth function whose sup-norm distance 
 from the  continuous function 
 $|s+3w|$ is less than $1/j$.) 
  Then Lemma \ref{nearly} tells us that   
 $$\int_M (s+3w)^2 d\mu = \inf_j \int_Mu_j^2 d\mu \geq 72\pi^2 (a^+)^2,$$
 as claimed.
   \end{proof}
   
   Note that  Lemma \ref{chiave} essentially reproves  
   \cite[Theorem 2.3]{lric}, but does so in a much more elementary manner. 
   What has been lost here, however, 
   is a precise understanding of what happens when the inequality is saturated. 
   For the purpose of finding obstructions to the existence  of
   Einstein metrics, however, this shortcoming will  in practice turn out to be  irrelevant.
   
\begin{prop}\label{goody}
  Let $M$ be a smooth compact oriented $4$-manifold with $b_+(M) \geq 2$. Then every 
   Riemannian metric $g$ on $M$ satisfies 
   $$\|s\| + \sqrt{6} \|W_+\|  \geq 6\sqrt{2}\pi\alpha (M),$$
   where $\|\cdot\|$ denotes the $L^2$ norm with respect to $g$.  
If equality occurs, moreover, then the two largest  eigenvalues of 
   $W_+: \Lambda^+\to \Lambda^+$
   are equal at each $x\in M$, and $|W_+|$ is a constant multiple of the scalar curvature $s$.
\end{prop} 
\begin{proof}
By Lemma \ref{chiave}, we have
$$\|s + 3w\|\geq  6\pi\sqrt{2 (a^+)^2},$$
so the triangle inequality tells us that 
$$
\|s\| + 3\| w\| \geq \|s + 3w\|\geq  6\pi\sqrt{2 (a^+)^2},
$$
and that equality can only hold if $w$ and $s$ are proportional, 
as vectors in $L^2$. 

Now if $\alpha (M)=0$, the claim is of course trivial. Otherwise, 
the definition of $\alpha (M)$ guarantees that
 for every metric $g$  there is a monopole class $a$ with 
$(a^+)^2 \geq \alpha^2 (M)$.  Now invoking our calculations for  this choice of $a$, we then have  
$$\|s\| + 3\| w\| \geq  6\pi\sqrt{2}\alpha (M).$$ 
Moreover, since a $4$-manifold $M$ with $\alpha (M) \neq 0$ 
cannot admit metrics with $s\equiv 0$, equality can only occur if
$w$ is a constant times $s$.  

On the other hand, because $W_+$ is a trace-free
endomorphism of $\Lambda^+$, we have the point-wise inequality 
 $$\sqrt{\frac{2}{3}}|W_+| \geq |w|,$$
 with equality only when the two largest eigenvalues  of $W_+$ are equal.
 Hence 
 $$\|s\| + \sqrt{6} \|W_+\|  \geq \|s\| + 3\| w\| \geq 6\pi\sqrt{2}\alpha (M).$$
 The resulting inequality is now strict unless the largest eigenvalue  of 
 $W_+$ has multiplicity $\geq 2$ everywhere, and unless
 the functions $|W_+|$ is a constant multiple of the scalar curvature $s$. 
\end{proof}

\begin{lem}\label{able}
Let $M$ be a smooth compact oriented $4$-manifold with $b_+(M) \geq 2$.
Then every Riemannian metric $g$ on $M$ satisfies
$$\frac{1}{4\pi^2}\int_M \left(\frac{s^2}{24} + 2|W_+|^2\right)d\mu \geq \frac{2}{3} \alpha^2 (M).$$
Moreover, equality could  only possibly hold for metrics for which $s\equiv  -8\sqrt{6}|W_+|$ and 
for which the  largest
eigenvalue of $W_+$ has multiplicity two   at each point of $M$. 
\end{lem}
\begin{proof}
By Proposition \ref{goody}, we have 
  $$(1 , {\frac{1}{\sqrt{8}}}) \cdot \left(\|s\|, 
 \sqrt{48} \|W_{+}\|\right) \geq 6\sqrt{2}\pi\alpha (M),$$
where the left-hand side is to be read as a dot product 
 in ${\mathbb R}^{2}$. Applying  the Cauchy-Schwarz equality to this dot product now gives 
us 
 $$ \left(1+\frac{1}{8}\right)^{1/2}
 \left(\|s\|^2+48\|W_{+}\|^{2} \right)^{1/2}
 \geq (72\pi^2 \alpha^2 (M))^{1/2},$$
 so that 
\begin{equation}
\label{loose}
\frac{1}{4\pi^2}\int_M \left(\frac{s^2}{24} + 2|W_+|^2\right)d\mu \geq \frac{2}{3} \alpha^2 (M),
\end{equation}
as claimed. If equality holds, we must  have $\left(\|s\|, 
 \sqrt{48} \|W_{+}\|\right) \propto (1 , {\frac{1}{\sqrt{8}}})$ in addition to the 
previously deduced conditions, and the value of the ratio of $|W_+|$ and $s$
then follows  from this, together with the fact that there cannot be any metrics with 
$s\geq 0$ when $\alpha (M) \neq 0$. \end{proof}

\begin{lem} \label{enable} 
Let $M$ be a smooth compact oriented $4$-manifold with $b_+(M) \geq 2$,
and suppose that $g$ is a Riemannian metric on $M$ with constant scalar curvature 
and harmonic
self-dual Weyl curvature:
$$(\delta W_+)_{bcd}=-\nabla^a(W_+)_{abcd} = 0.$$
Then either $g$  satisfies  the strict inequality
$$
\frac{1}{4\pi^2}\int_M \left(\frac{s^2}{24} + 2|W_+|^2\right)d\mu > \frac{2}{3} \alpha^2 (M),
$$
or else $\alpha (M) =0$, $g$ is a hyper-K\"ahler metric, and $M$ is diffeomorphic to either
$K3$ or $T^4$. 
\end{lem}
\begin{proof}
Suppose that $g$ is a metric with constant $s$  
 which saturates  inequality (\ref{loose}).
 Then $|W_+|\equiv  -s/8\sqrt{6}$ is constant, and 
the  largest two 
eigenvalues of $W_+$ coincide   at each point of $M$. 
We now want to ask what this tells us if we also assume that $\delta W_+=0$. 

One way this may certainly happen is for $s$ and $W_+$ to both vanish identically. 
In this case, however, the Weitzenb\"ock formula (\ref{friend}) then 
implies that every self-dual harmonic form on $M$ is parallel,
and, since $b_+(M)\geq 2$ by assumption, it then follows that the orientable 
rank-$3$ vector bundle 
$\Lambda^+$ is trivialized by parallel sections, so that 
$g$ is hyper-K\"ahler. The classification of complex surfaces then 
tells us  that $M$ must  be diffeomorphic to $K3$ or $T^4$. 

On the other hand, the  $s\neq 0$ case  is ruled out by  an 
observation due to Derdzi{\'n}ski \cite{derd}. Indeed, in this case we would know that 
$W_+$ was a non-zero constant times $\omega \otimes
 \omega -\frac{1}{3}{\mathbf 1}_{\Lambda^+}$, where $\omega$ is a (sign-ambiguous)
self-dual $2$-form of norm $\sqrt{2}$ spanning the negative eigenspace of 
$W_+$. The harmonicity of $W_+$ thus implies that 
$$\omega_{ab}\nabla^a \omega_{cd} + \omega_{cd}\nabla^a \omega_{ab} =0.$$
But since $\omega$ has constant length, contraction with $\omega^{cd}$
tells us that $\delta\omega =0$, and plugging this back into the original equation 
then tells us that $\nabla \omega =0$. This shows that $(M,g)$ is locally K\"ahler.
 However,
any  K\"ahler manifold of real dimension $4$ satisfies $|s|\equiv 2\sqrt{6}|W_+|$, whereas
in our case  we can  already take it as given that 
$|s|\equiv 8\sqrt{6}|W_+|$. This contradiction eliminates the $s\neq 0$ case, and 
we are done. \end{proof}

This  now allows us to recover  \cite[Theorem 3.3]{lric}:

\begin{thm}\label{righto}
Let $M$ be a smooth compact oriented $4$-manifold with $b_+(M) \geq 2$.
If $M$ admits an  Einstein metric $g$, then 
$$
(2\chi + 3\tau ) (M) \geq \frac{2}{3} \alpha^2 (M),
$$
with equality only if both sides vanish, in which case 
$M$ must be  diffeomorphic to $K3$ or $T^4$, and $g$
must be  a hyper-K\"ahler metric. 
\end{thm}
\begin{proof}
If $(M,g)$ is a compact oriented $4$-dimensional Einstein manifold, 
the generalized Gauss-Bonnet formula  (\ref{gb}) tells us 
that the left-hand side of (\ref{loose}) is given by $(2\chi+3\tau )(M)$, and the 
desired inequality therefore follows. However, 
 the Bianchi identities also tell us that a $4$-dimensional 
Einstein manifold has  constant $s$ and harmonic $W_+$, so Lemma
\ref{enable} tells us that this inequality can only 
be saturated when $g$ is hyper-K\"ahler. 
\end{proof}

\section{An Illustration}

Theorem \ref{righto}  gives us an 
obstruction to the existence of Einstein metrics 
involving the differential topology,
rather than just the homeomorphism type,
of the smooth manifold in question. 
The following example should
help to clarify this point.

Let $X$ be a triple cyclic cover of of $\CP_2$, ramified at a 
non-singular complex curve $B$ of degree $6$.
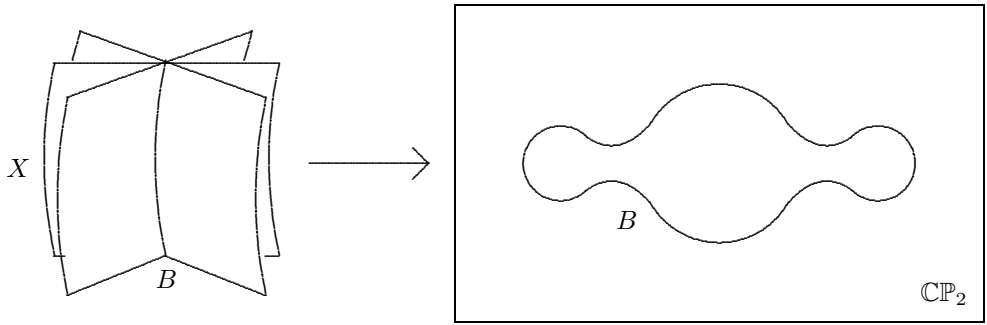
\begin{center}
\mbox{
\beginpicture
\setplotarea x from 30 to 400, y from -30 to 120
{\setlinear 
\plot 55 88 140 88 /
\plot  55 15 59 15 /
\plot  140 15 135 15 /
\plot 65 100 135 75 /
\plot 60 75  130 100 /
\plot 60 0 97 15 /
\plot 97 15 135 0 /
\plot 130 100 127 89 /
\plot 65 100 62 89 /
}
{\setquadratic
\plot 94 89  90  52  94 15 /
\plot 132 75 128 37 132 0 /
\plot 57 75 53 37 57 0 /
\plot 52 88 48 52 52 15 /
\plot 137 88 133 52 137 15 /
}
\putrectangle corners at 200 110 and 400 -10
\circulararc 270 degrees from 250 60  center at 240 50
\circulararc -270 degrees from 350 60 center at 360 50
\circulararc 120 degrees from 326 65 center at 300 50 
\circulararc -120 degrees from 326 35  center at 300 50 
\put {${\mathbb C \mathbb P}_2$} [B1] at 385 -2
\put {$B$} [B1] at 265 25
\put {$B$} [B1] at 91 3
\put {$X$} [B1] at 35 45
\arrow <6pt> [1,2] from 145 50  to 190 50
{\setquadratic
\plot 350 60 338 57 326 65 /
\plot 350 40 338 43 326 35 /
\plot 250 40 262 43 274 35 /
\plot 250 60 262 57 274 65 /
}
\endpicture
}
\end{center}
To be more explicit, we could, for example, 
 take  $X$ to consist of those
elements $\eta$ of the ${\mathcal O} (2)$ line bundle over $\CP_2=\{ [x:y:z]\}$ which satisfy the 
equation 
$$
\eta^3 = x^6 + y^6 + z^6 , 
$$
where the right-hand side is of course interpreted as a section of ${\mathcal O}(6)$. 
The canonical line bundle of the compact complex surface
 $X$ is then exactly the pull-back of the ${\mathcal O} (1)$ line bundle on $\CP_2$,
and from this we read off that 
$c_1^2 (X) = 3 \cdot 1^2= 3$, and $h^{2,0}(X) = h^0 (\CP_2 , {\mathcal O} (1))= 3$. 
Moreover, an easy application of the Lefschetz hyperplane-section theorem 
 tells us that $X$ is  simply connected.

Now  
let $M$ be obtained from $X$ by blowing up a point, so that 
$$M \approx X \# \overline{\CP}_2 .$$
Then $\alpha^2 (M) = c_1^2 (X) = 3$, whereas 
$(2\chi + 3\tau ) (M) = c_1^2 (M) = c_1^2 (X) - 1 = 2$. 
Since $(2\chi + 3\tau ) (M) = \frac{2}{3}\alpha^2(M)\neq 0$, Theorem \ref{righto}
therefore 
tells us that  the smooth compact 
$4$-manifold $M$  cannot admit an Einstein metric.

Next, for the sake of  comparison, let $N$ be the double branched cover of $\CP_2$, ramified at a 
non-singular complex curve $B^\prime$
of degree $8$. 
To be more explicit,   one could take  $N$ to consist of those
elements $\zeta$ of the ${\mathcal O} (4)$ line bundle over $\CP_2$ which satisfy the 
equation 
$$
\zeta^2 = x^8 + y^8 + z^8 .
$$
Once again, this branched cover of $\CP_2$ is simply connected.
Careful inspection also reveals that, in this example,  
 the canonical 
line bundle  is once again exactly the pull-back of the ${\mathcal O} (1)$ line bundle on $\CP_2$;
thus we  read off  that 
$c_1^2 (N) = 2 \cdot 1^2= 2$, and $h^{2,0}(N) = h^0 (\CP_2 , {\mathcal O} (1))= 3$,
which is to say that these two numerical invariants are exactly the same for 
$M$ and $N$. 
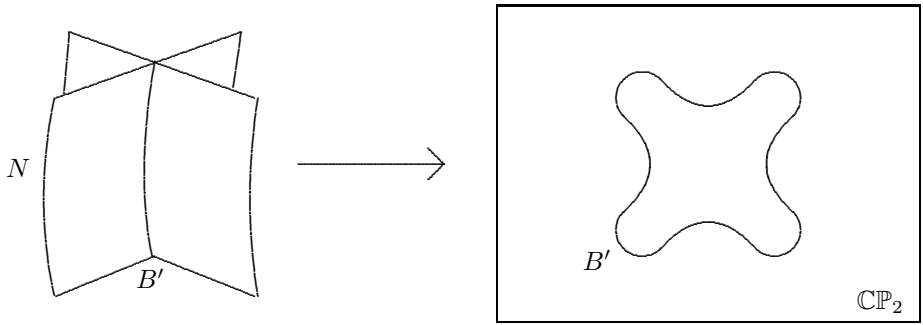
\begin{center}
\mbox{
\beginpicture
\setplotarea x from 30 to 400, y from -30 to 120
{\setlinear 
\plot 65 100 135 75 /
\plot 60 75  130 100 /
\plot 60 0 97 15 /
\plot 97 15 135 0 /
\plot 130 100 127 79 /
\plot 65 100 63 77 /
}
{\setquadratic
\plot 133 75 130 37 133 0 /
\plot 56 75 52 37 56 0 /
\plot 94 89 90 49 93 15 /
}
\putrectangle corners at 220 110 and 380 -10
\circulararc 180 degrees from 332 68  center at 325 75
\circulararc -180 degrees from 332 32  center at 325 25
\circulararc -180 degrees from 268 68  center at 275 75
\circulararc 180 degrees from 268 32  center at 275 25
\put {${\mathbb C \mathbb P}_2$} [B1] at 365 -5
\put {$B^\prime$} [B1] at 258 10
\put {$B^\prime$} [B1] at 89 3
\arrow <6pt> [1,2] from 145 50  to 200 50
\put {$N$} [B1] at 39 45
{\setquadratic
\plot 282 82  300 72  318 82 /
\plot 282 18  300 28  318 18 /
\plot  332 32  322 50   332  68 /
\plot  268 32  278 50   268  68 /
}
\endpicture
}
\end{center}
But in another respect, $N$ is wildly unlike  from $M$. Indeed,  because it
is a  compact complex surface with $c_1 < 0$.
Theorem \ref{ay} tells us that $N$, unlike $M$, admits  an Einstein  metric.

However, each  of these compact complex surfaces
is simply connected and non-spin. Moreover,  both of them 
have $c_1^2 = 2$ and $h^{2,0}=3$, and from this it can be deduced that
both have  
$b_+=7$ and $b_-=37$. Thus Theorem \ref{fdmn}  
tells us that both $M$ and $N$ are homeomorphic to 
$$7\CP_2 \# 37\overline{\CP}_2.$$
This topological $4$-manifold therefore admits one
smooth structure for which there is an Einstein metric, and 
yet another smooth structure for which no
such metric can exist.

\section{The Big Picture} 
\label{entropy}

The problem of understanding  Einstein metrics
on compact $4$-manifolds has been an extremely active area of 
research in recent years, and the Seiberg-Witten techniques which
have been highlighted in  this article  do not by any means 
represent the only strand of thought which has led 
to significant progress. It is therefore appropriate that this
article conclude with a rough indication of some of these other developments, 
without pretending to offer a definitive survey of the subject. 

One such important strand of thought grew  out of the work of Gromov \cite{grom},
who, for example, first introduced  
minimal volume invariants of  a smooth manifold.  If $M$ is a smooth compact
$n$-manifold, one such invariant is 
$$\Vol_r (M) =\inf \{ \Vol (M,g )~|~ r_g \geq -(n-1)g \},$$
and 
is sometimes called the {\em minimal volume with respect to Ricci curvature}. 
 Gromov showed that  $\Vol_r (M)$ is positive for
 certain compact manifolds with infinite fundamental 
 group by giving a lower bound for it in terms of  a 
 homotopy invariant which he called the {\em simplicial
 volume}. He then went on to observe that this implies that 
 there are obstructions to the 
existence of Einstein metrics in dimension $4$ which are not detected
 by the Hitchin-Thorpe inequality. 

 While Gromov's bounds are actually quite weak in practice,
 they nonetheless represented a breakthrough
in the subject, and eventually  led to a  hunt for
sharper estimates of the same flavor. This culminated in the work of  
Besson, Courtois, and Gallot \cite{bcg}, who related  Gromov's work to 
the entropy of the geodesics flow. One of their most remarkable
results goes as follows:

\begin{thm}[Besson-Courtois-Gallot]  \label{qdeg} 
Let $(X,g_0)$ be a compact hyperbolic $n$-manifold, $n > 2$, 
and let $M$ be a compact
manifold of the same dimension. Let $f: M\to X$ be any 
smooth map. Then \begin{equation}
\label{minvol}
\Vol_r (M) \geq \deg (f)~ \Vol (X,g_0),
\end{equation}
where $\deg (f)$ denotes the degree of $f$. 
Moreover, if equality holds,  and if the infimum (\ref{minvol})
is achieved by some metric $g$, then
$(M,g)$ is an isometric  Riemannian covering of $(X,g_0)$, 
with covering map $M\to X$  homotopic to $f$. 
\end{thm}

When $n=4$
and $f$ is the identity map $X\to X$, 
this  implies  that $g_0$ is the only Einstein metric on $X$,
up to  rescalings and diffeomorphisms.
Moreover,  this same result also gives rise to new obstructions
to the existence of Einstein metrics. 
Indeed, 
Sambusetti \cite{samba} used this method to prove the following:

\begin{thm}[Sambusetti] \label{samba} 
Any integer pair $(\chi,\tau)$ with $\chi\equiv  \tau\bmod 2$
can be realized as the Euler characteristic   and signature
of a smooth   compact  oriented $4$-manifold $M$  
which does not admit any 
 Einstein metrics.
\end{thm}

This  nicely highlights how much there is to be 
said about the subject beyond 
the Hitchin-Thorpe inequality.
Note, however, that 
Sambusetti's examples all have huge fundamental group, and 
are   never even {\em homotopy equivalent} to an Einstein manifold.

Anderson's work \cite{anderson}
on Gromov-Hausdorff limits with Ricci-curvature bounds 
represents another major area of progress in understanding $4$-dimensional
Einstein manifolds.   Anderson shows that any sequence of unit-volume
Einstein manifolds of bounded Euler characteristic and
{\em bounded diameter} 
must necessarily have a subsequence which converges in the Gromov-Hausdorff 
sense to 
a compact Einstein orbifold. In particular, this implies that 
 a topological 
$4$-manifold can only admit finitely many  smooth structures for which there exists 
an Einstein metric of unit volume and diameter $<D$. 

It is interesting to reconsider the  K\"ahler-Einstein case in this light. 
Mumford's finiteness theorem \cite{bpv} implies that a $4$-manifold
can only admit finitely many smooth structures for which there exist
{\em K\"ahler}-Einstein metrics, although 
this finite number
can be arbitrarily large  \cite{kot,salvetti}.
Notice that this conclusion holds without the imposition of  
an extraneous 
diameter bound. On the other hand,  inspection of the K\"ahler-Einstein case also reveals  that 
 the diameter can tend to infinity even for sequences of unit-volume 
 Einstein metrics  on a fixed smooth $4$-manifold, so a Mumford-type  finiteness statement certainly 
 cannot  be deduced  from a
 compactness result like Anderson's.

 While we do not know at present whether there can only be finitely many
 smooth structures admitting  Einstein metrics on a fixed topological $4$-manifold,
 we {\em do} at least know that there are often infinitely many smooth structures
 for which Einstein metrics {\em don't} exist, even for simply connected 
 manifolds for which the Hitchin-Thorpe inequality is strict \cite{il1}; cf. \cite{kot}.
In the non-spin case, such an assertion can actually be made for  a
large sector of  choices for $(\chi , \tau )$.
 The spin case is much more delicate, however, 
and the range of homeotypes for which one can make such an assertion 
is, at present,  much more tightly constrained.

%

%






%


\begin{thebibliography}{10}

\bibitem{anderson}
M.T. Anderson.
\newblock The ${L}^{2}$ structure of moduli spaces of {E}instein metrics on
  $4$-manifolds.
\newblock {\em Geom. Func. An.}, 2:29--89, 1992.

\bibitem{mta}
M.T. Anderson.
\newblock Dehn filling and {E}instein metrics in higher dimensions.
\newblock e-print math.DG/0303260, 2003.

\bibitem{aubin}
T. Aubin.
\newblock Equations du type {M}onge-{A}mp\`{e}re sur les vari\'{e}t\'{e}s
  {K\"a}hl{\'e}riennes compactes.
\newblock {\em C. R. Acad. Sci. Paris}, 283A:119--121, 1976.

\bibitem{bpv}
W.~Barth, C.~Peters, and A.~Van de~Ven.
\newblock {\em Compact Complex Surfaces}.
\newblock Springer-Verlag, 1984.

\bibitem{baufu}
S.~Bauer and M.~Furuta.
\newblock A stable cohomotopy refinement of {S}eiberg-{W}itten invariants: I.
\newblock Inv. Math. to appear.

\bibitem{bauer2}
S. Bauer.
\newblock A stable cohomotopy refinement of {S}eiberg-{W}itten invariants:
  {II}.
\newblock Inv. Math. to appear.

\bibitem{bes}
A.~Besse.
\newblock {\em {E}instein Manifolds}.
\newblock Springer-Verlag, 1987.

\bibitem{bcg}
G.~Besson, G.~Courtois, and S.~Gallot.
\newblock Entropies et rigidit{\'e}s des espaces localement sym{\'e}triques de
  courbure strictement n{\'e}gative.
\newblock {\em Geom. and Func. An.}, 5:731--799, 1995.

\bibitem{borel}
A. Borel.
\newblock Compact {C}lifford-{K}lein forms of symmetric spaces.
\newblock {\em Topology}, 2:111--122, 1963.

\bibitem{bourg}
J.-P. Bourguignon.
\newblock Les vari\'et\'es de dimension $4$\ \`a signature non nulle dont la
  courbure est harmonique sont d'{E}instein.
\newblock {\em Invent. Math.}, 63(2):263--286, 1981.

\bibitem{derd}
A. Derdzi{\'n}ski.
\newblock Self-dual {K}\"ahler manifolds and {E}instein manifolds of dimension
  four.
\newblock {\em Compositio Math.}, 49(3):405--433, 1983.

\bibitem{donaldson}
S.~K. Donaldson.
\newblock An application of gauge theory to four-dimensional topology.
\newblock {\em J. Differential Geom.}, 18:279--315, 1983.

\bibitem{don}
S.~K. Donaldson and P.~B. Kronheimer.
\newblock {\em The Geometry of Four-Manifolds}.
\newblock Oxford University Press, Oxford, 1990.

\bibitem{freedman}
M.~Freedman.
\newblock On the topology of 4-manifolds.
\newblock {\em J. Differential Geom.}, 17:357--454, 1982.

\bibitem{fur108}
M.~Furuta.
\newblock Monopole equation and the {$\frac{11}8$}-conjecture.
\newblock {\em Math. Res. Lett.}, 8(3):279--291, 2001.

\bibitem{gost}
R.E. Gompf and A.I. Stipsicz.
\newblock {\em $4$-manifolds and {K}irby calculus}.
\newblock American Mathematical Society, Providence, RI, 1999.

\bibitem{grom}
M.~Gromov.
\newblock Volume and bounded cohomology.
\newblock {\em Publ. Math. IHES}, 56:5--99, 1982.

\bibitem{gvln}
M.~Gromov and H.B. Lawson.
\newblock The classification of simply connected manifolds of positive scalar
  curvature.
\newblock {\em Ann. Math.}, 111:423--434, 1980.

\bibitem{hitharm}
N.J. Hitchin.
\newblock Harmonic spinors.
\newblock {\em Advances in Mathematics}, 14:1--55, 1974.

\bibitem{hit}
N.J. Hitchin.
\newblock On compact four-dimensional {E}instein manifolds.
\newblock {\em J. Differential Geom.}, 9:435--442, 1974.

\bibitem{il2}
M.~Ishida and C.~LeBrun.
\newblock Curvature, connected sums, and {S}eiberg-{W}itten theory.
\newblock to appear in {C}omm. {A}nal. {G}eom.

\bibitem{il1}
M.~Ishida and C.~LeBrun.
\newblock Spin manifolds, {E}instein metrics, and differential topology.
\newblock {\em Math. Res. Lett.}, 9:229--240, 2002.

\bibitem{kot}
D.~Kotschick.
\newblock {E}instein metrics and smooth structures.
\newblock {\em Geom. Topol.}, 2:1--10, 1998.

\bibitem{K}
P.B.~Kronheimer.
\newblock Minimal genus in ${S}\sp 1\times {M}\sp 3$.
\newblock {\em Invent. Math.}, 135(1):45--61, 1999.

\bibitem{lawmic}
H.B. Lawson and M.L. Michelsohn.
\newblock {\em Spin Geometry}.
\newblock Princeton University Press, 1989.

\bibitem{lpm}
C.~LeBrun.
\newblock Polarized 4-manifolds, extremal {K\"a}hler metrics, and
  {S}eiberg-{W}itten theory.
\newblock {\em Math. Res. Lett.}, 2:653--662, 1995.

\bibitem{lno}
C.~LeBrun.
\newblock Four-manifolds without {E}instein metrics.
\newblock {\em Math. Res. Lett.}, 3:133--147, 1996.

\bibitem{lric}
C.~LeBrun.
\newblock Ricci curvature, minimal volumes, and {S}eiberg-{W}itten theory.
\newblock {\em Inv. Math.}, 145:279--316, 2001.

\bibitem{lic}
A.~Lichnerowicz.
\newblock Spineurs harmoniques.
\newblock {\em C.R. Acad. Sci. Paris}, 257:7--9, 1963.

\bibitem{ozsz}
P. Ozsv{\'a}th and Z. Szab{\'o}.
\newblock Higher type adjunction inequalities in {S}eiberg-{W}itten theory.
\newblock {\em J. Differential Geom.}, 55(3):385--440, 2000.

\bibitem{page}
D.~Page.
\newblock A compact rotating gravitational instanton.
\newblock {\em Phys. Lett.}, 79B:235--238, 1979.

\bibitem{pr2}
R. Penrose and W. Rindler.
\newblock {\em Spinors and space-time. {V}ol. 2}.
\newblock Cambridge University Press, Cambridge, 1986.
\newblock Spinor and twistor methods in space-time geometry.

\bibitem{salvetti}
M. Salvetti.
\newblock On the number of nonequivalent differentiable structures on
  $4$-manifolds.
\newblock {\em Manuscripta Math.}, 63(2):157--171, 1989.

\bibitem{samba}
A.~Sambusetti.
\newblock An obstruction to the existence of {E}instein metrics on 4-manifolds.
\newblock {\em C. R. Acad. Sci. Paris}, 322:1213--1218, 1996.

\bibitem{taubes}
C.H. Taubes.
\newblock The {S}eiberg-{W}itten invariants and symplectic forms.
\newblock {\em Math. Res. Lett.}, 1:809--822, 1994.

\bibitem{taubes3}
C.H. Taubes.
\newblock The {S}eiberg-{W}itten and {G}romov invariants.
\newblock {\em Math. Res. Lett.}, 2:221--238, 1995.

\bibitem{tho}
J.A. Thorpe.
\newblock Some remarks on the {G}auss-{B}onnet formula.
\newblock {\em J. Math. Mech.}, 18:779--786, 1969.

\bibitem{tian}
G.~Tian.
\newblock On {C}alabi's conjecture for complex surfaces with positive first
  chern class.
\newblock {\em Inv. Math.}, 101:101--172, 1990.

\bibitem{witten}
E.~Witten.
\newblock Monopoles and four-manifolds.
\newblock {\em Math. Res. Lett.}, 1:809--822, 1994.

\bibitem{yau}
S.-T. Yau.
\newblock {C}alabi's conjecture and some new results in algebraic geometry.
\newblock {\em Proc. Nat. Acad. USA}, 74:1789--1799, 1977.

\end{thebibliography}
\end{document}